\documentclass[11pt]{amsart}
\usepackage{amsmath,amssymb,amsthm}
\usepackage[latin1]{inputenc}
\usepackage{version,tabularx,multicol}
\usepackage{graphicx,float}

\newcommand{\reff}[1]{(\ref{#1})}

\theoremstyle{plain}
\newtheorem{theo}{Theorem}[section]

\newtheorem{cor}[theo]{Corollary}
\newtheorem{prop}[theo]{Proposition}

\newtheorem{lem}[theo]{Lemma}
\newtheorem{defi}[theo]{Definition}
\theoremstyle{remark}
\newtheorem{rem}[theo]{Remark}

\newcommand{\hm}{{H_{\text{max}}}}
\newcommand{\htl}{{H_{[T_m,L_m]}}}
\newcommand{\tm}{{T_{\text{max}}}}
\newcommand{\Em}{{\eta_{\tm}}}
\newcommand{\Rm}{{\rho_{\tm}}}
\newcommand{\km}{{\kappa_{\text{max}}}}
\newcommand{\pe}{\pi_{\text{Eve}}}
\newcommand{\pt}{\pi}
\newcommand{\Pt}{\psi}
\newcommand{\Pe}{\psi_{\text{Eve}}}
\newcommand{\at}{\alpha_0}
\renewcommand{\ae}{\alpha_{\text{Eve}}}
\newcommand{\aimm}{\alpha_{\text{Imm}}}
\newcommand{\mue}{\mu_{\text{Eve}}}
\newcommand{\mut}{\mu_{\text{Total}}}

\newcommand{\cb}{{\mathcal B}}

\newcommand{\ce}{{\mathcal E}}

\newcommand{\cn}{{\mathcal N}}
\newcommand{\cm}{{\mathcal M}}

\newcommand{\E}{{\mathbb E}}

\newcommand{\N}{{\mathbb N}}
\renewcommand{\P}{{\mathbb P}}

\newcommand{\Q}{{\mathbb Q}}
\newcommand{\R}{{\mathbb R}}

\newcommand{\rN}{{\rm N}}
\newcommand{\rP}{{\rm P}}
\newcommand{\rE}{{\rm E}}

\newcommand{\ind}{{\bf 1}}

\newcommand{\Supp}{{\rm Supp}\;}

\newcommand{\expp}[1]{\mathop {\mathrm{e}^{ #1}}}

\begin{document}

\title[Williams' decomposition]{Williams' decomposition of the
  Lévy continuum random tree and 
  simultaneous extinction probability for populations with
 neutral mutations}

\date{\today}
\author{Romain Abraham}

\address{
Romain Abraham,
MAPMO, F\'ed\'eration Denis Poisson, Université d'Orléans,
B.P. 6759,
45067 Orléans cedex 2,
France.
}
  
\email{romain.abraham@univ-orleans.fr}

\author{Jean-François Delmas}

\address{
Jean-Fran\c cois Delmas,
CERMICS, Univ. Paris-Est,  6-8
av. Blaise Pascal, 
  Champs-sur-Marne, 77455 Marne La Vallée, France.}

\email{delmas@cermics.enpc.fr}

\begin{abstract}
  We  consider an  initial Eve-population  and a  population  of neutral
  mutants,  such that  the total  population  dies  out  in finite  time.  We
  describe the evolution of  the Eve-population and the total population
  with continuous  state branching  processes, and the  neutral mutation
  procedure  can  be  seen  as  an immigration  process  with  intensity
  proportional  to the  size of  the  population. First  we establish  a
  Williams' decomposition of the  genealogy of the total population given
  by a continuum random tree, according to the ancestral lineage of the
  last individual  alive. This allows us to give a closed  formula for the
  probability of  simultaneous extinction of the  Eve-population and the
  total population.
\end{abstract}



\maketitle

\section{Introduction}

We consider an initial Eve-population whose size evolves as a continuous
state  branching process  (CB), $Y^0=(Y^0_t,  t\geq 0)$,  with branching
mechanism $\Pe$.  We assume this  population gives birth to a population
of irreversible  mutants. The new mutants  population can be  seen as an
immigration  process   with  rate  proportional  to  the   size  of  the
Eve-population.   We assume  the  mutations are  neutral,  so that  this
second population  evolves according to the same  branching mechanism as
the Eve-population.   This population of  mutants gives birth also  to a
population of other irreversible  mutants, with rate proportional to its
size,  and   so  on.   In  \cite{ad:cbpimcsbpi},  we   proved  that  the
distribution of the total population size $Y=(Y_t, t\geq 0)$, which is a
CB with immigration (CBI) proportional to its own size, is in fact a CB,
whose branching mechanism $\Pt$ depends on the immigration intensity.  
The joint law of $(Y^0,Y)$ is characterized by its Laplace transform,
see Section \ref{sec:Lap-trans}. This model can also be viewed as a
special case of multitype CB, with two types 0 and 1, the individuals
of type 0 giving birth to offsprings of type 0 or
1, whereas individuals of type 1 only have type 1 offsprings, see \cite{m:ipsmgwt, b:saptpgwpnm} for recent related works.

In the  particular case of $Y$  being a sub-critical or  critical CB with
quadratic  branching mechanism  ($\Pt(u)=\at u+  \beta  u^2$, $\beta>0$,
$\at\geq 0$), the probability for the Eve-population to disappear at the
same time as the whole population is known, see \cite{w:bprtsbm} for the
critical  case, $\at=0$, or  Section 5  in \cite{ad:cbpimcsbpi}  for the
sub-critical case, $\at> 0$. Our aim is to extend those results for
 the large
class of CB with unbounded total variation and  a.s. extinction.  Formulas   given  in
\cite{ad:cbpimcsbpi} could certainly be extended to a general branching
mechanism, but first computations seem to be rather involved.

In fact,  to compute  those quantities,  we choose here to  rely   on the
description of the  genealogy of sub-critical or critical  CB introduced  by Le  Gall  and  Le  Jan
\cite{lglj:bplpep}  and   developed  later  by  Duquesne   and  Le  Gall
\cite{dlg:rtlpsbp}, see also Lambert \cite{l:gcsbpi} for the genealogy
of CBI with constant immigration rate. Le Gall and Le Jan defined via a
L\'evy process $X$ the so-called height process $H=(H_t,t\ge 0)$ which
codes a continuum random tree (CRT) that describes the genealogy of the
CB (see the next section for the definition of $H$ and the coding of the
CRT).
Initially, the CRT was introduced by Aldous \cite{a:crt1} in the quadratic case:
$\psi(\lambda)=\lambda^2$. Except in this quadratic case, the height
 process $H$ is not Markov and so is
difficult to handle. That is why they also introduce a measure-valued
Markov process $(\rho_t,t\ge 0)$ called the exploration process and
such that the closed support of the measure $\rho_t$ is $[0,H_t]$ (see
also the next section for the definition of the exploration process). 

We shall  be interested in the case  where a.s.  the
extinction of the whole population  holds in finite time.  The branching
mechanism of the total population, $Y$, is given by: for $\lambda\geq 0$,
\begin{equation}
   \label{eq:def_psi}
\Pt(\lambda)=\at\lambda+\beta \lambda^2+ \int_{(0,\infty)}\pt(d\ell)
\left(\expp{-\lambda\ell}-1+\lambda \ell\right),  
\end{equation}
where $\at\geq  0$, $\beta\geq 0$ and  $\pt$ is a  Radon measure on
$(0,\infty  )$ such  that $\int_{(0,\infty  )} (\ell  \wedge  \ell^2) \;
\pt(d\ell)<\infty $. We shall assume  that $Y$ is of infinite variation,
that   is   $\beta>0$  or   $\int_{(0,1)}   \ell  \pt(d\ell)=\infty   $.
   We shall
assume that  a.s.  the extinction of  $Y$ in finite time holds, that is,
see Corollary 1.4.2 in \cite{dlg:rtlpsbp},
we assume that 
\begin{equation}
   \label{eq:contH}\int^{+\infty }  \frac{dv}{\Pt(v)}<\infty .
\end{equation}  
We suppose that the process $Y$ is the canonical process on the
Skorokhod space $\mathbb{D}(\R_+,\R_+)$ of c\`adl\`ag paths and that
the pair $(Y,Y^0)$ is the canonical process on the space $ \mathbb{D}(\R_+,\R_+)^2$.
Let
$\rP_x$ denote  the law of the pair $(Y,Y^0)$ (see
\cite{ad:cbpimcsbpi}) started at  $(Y_0,Y^0_0)=(x,x)$. The probability
measure $\rP_x$ is infinitely divisible and hence admits a canonical
measure $\rN$: it is a $\sigma$-finite measure on $
\mathbb{D}(\R_+,\R_+)^2$ such that
$$(Y,Y^0)\overset{(d)}{=}\sum_{i\in I}(Y^i,Y^{0,i})$$
where $((Y^i,Y^{0,i}), i\in I)$ are the atoms of a Poisson measure on $
\mathbb{D}(\R_+,\R_+)^2$ with intensity $x\rN(dY,dY^0)$. In particular,
we have 
\begin{equation}
   \label{eq:br-Y}
\rE_x[\expp{-\lambda Y_t}]=\exp (- x \rN[1-\expp{-\lambda Y_t}])
\end{equation}
 and
$u(\lambda,t)= \rN[1-\expp{-\lambda Y_t}]$ is the unique non-negative
solution of 
\begin{equation}
   \label{eq:u-lt}
\int_{u(\lambda, t)}^\lambda \frac{dv}{\Pt(v)}=t, \quad\text{for $t\geq 0$
  and $\lambda\geq 0$.}
\end{equation}
Let $\tau_Y=\inf\{t>0; Y_t=0\}$ be the extinction time of $Y$.  Letting
$\lambda $ go to $\infty $ in the previous equalities leads to 
\[
\rP_x(\tau_Y<t)=\exp -x \rN[\tau_Y\geq t],
\]
where the positive function $c(t)=\rN[\tau_Y\geq t]$ solves
\begin{equation}
   \label{eq:def-c}
 \int_{c(t)}^\infty
\frac{dv}{\Pt(v)}=t, \quad\text{for $t>0$.}
\end{equation}

Let us consider the exploration process $(\rho_t,t\ge 0)$ associated
with this CB. We denote by $\N$ its excursion measure. Recall that the
closed support of the measure $\rho_t$ is $[0,H_t]$, where $H$ is the
height process. Let $L^a$ be the
total local time at level $a$ of the height process $H$ (well-defined
under $\N$). Then, the process $(L^a,a\ge 0)$ under $\N$ has the same
distribution as the CB $Y$ under $\rN$.

We decompose the exploration process, under the excursion measure,
according to the maximum of the height process. In terms of the CRT,
this means that we consider the
longest rooted  branch of the  CRT and describe how the different
subtrees are grafted along that branch, see Theorem  \ref{theo:rmax}.
When the branching mechanism is quadratic, the height process $ H$  is a
Brownian excursion and  the exploration process
$\rho_t$ is,  up to  a constant, the  Lebesgue measure on  $[0,H_t]$. In
that  case,   this  decomposition  corresponds   to  
Williams' original decomposition of the Brownian excursion (see
\cite{w:pdcltodd}).
This kind of tree decomposition with respect to a particular branch
(or a particular subtree) is not new, let us cite \cite{gk:slgwtpiv,
  pw:gbf} for instance, or \cite{s:omsbmcne, sv:cemsbm, ew:dscs} for related
works on superprocesses.

We present in the introduction a Poisson decomposition for the CB
only, and
we refer to Theorem  \ref{theo:rmax} for the decomposition of the
exploration process. 
Conditionally on the extinction time $\tau_Y$ equal to $m$, we can
represent the process $Y$ as the sum of the descendants of the ancestors
of the last individual alive. More precisely, let 
$\displaystyle \cn'(d\ell,dt)= \sum_{i\in I} \delta_{(\ell_i, 
  t_i)}(d\ell,dt) $ be a Poisson point measure with intensity 
\[
\ind_{[0,m)}(t) \expp{-\ell c(m-t)}  \ell \pt(d\ell)dt,
\]
and
\begin{equation}
   \label{eq:def-km}
\km(dt)=\sum_{i\in I} \ell_i
\delta_{t_i}(dt) + 2\beta \ind_{[0,m)} (t) dt.
\end{equation}

Let $\rN_t(dY)$ denote the law of
$(Y(s-t), s\ge t)$ under $\rN$ and 
$\displaystyle  \sum_{j\in J} \delta_{(t_j, Y^j)}$ be, conditionally on
 $\cn'$, a Poisson point
measure with intensity 
\[
\km(dt) \rN_t[dY, \ind_{\{\tau_Y\le m\}}]
\]
where $\rN_t[dY, \ind_{\{\tau_Y\le m\}}]$ denotes the restriction of the
measure $N_t$ to the event $\{\tau_Y\le m\}$.

The next result is  a direct consequence of Theorem
\ref{theo:rmax}.
\begin{prop}
\label{prop:Ylaw}
  The process  $ \sum_{j\in J} Y^j$  is distributed as  $Y$ under $\rN$,
  conditionally on $\{\tau_Y=m\}$.
\end{prop}
Let  $\tau_{Y^0}=\inf\{t>0; Y^0_t=0\}$ be the extinction time of the
Eve-population. 
In the particular case where the branching mechanism of the 
Eve-population is given by a shift of $\Pt$:
\begin{equation}
\label{eq:shift}
\Pe(\cdot)=\Pt(\theta+\cdot)-\Pt(\theta),
\end{equation}
for  some $\theta>0$  and $\beta=0$,  the  pruning procedure
developed in  \cite{ad:falp} gives that the nodes of width $\ell_i$
correspond to a mutation with probability $1- \expp{-\theta \ell_i}$. As
 $\beta=0$ there is no mutation on the skeleton of the CRT outside the
 nodes. In particular, we see simultaneous extinction of the whole
 population and the Eve-population if there is no mutation on the nodes
 in the ancestral lineage of the last individual alive. This happens,
 conditionally on $\km$, 
 with probability 
\[
\expp{-\theta \sum_{i\in I} \ell_i}. 
\]
Integrating w.r.t. the  law of $\cn$ gives that  the
probability of simultaneous extinction, conditionally on
$\{\tau_Y=m\}$, is under $\rN$, 
given by 
\begin{align*}
\rN[\tau_{Y^0}=m|\tau_Y=m]
&=\exp {-\int\ind_{[0,m)}(t) \expp{-\ell c(m-t)}  \ell \pt(d\ell)dt
  \left[1-\expp{-\theta \ell}\right]}\\
&= \exp{-\int_0^m
[  \Pt'(c(m-t)+\theta) - \Pt'(c(m-t)) ]\;dt }\\
&= \exp{-\int_0^m
  \phi'(c(t)) \;dt },
\end{align*}
where $\phi=\Pe-\Pt$. Now, using that the distribution of $(Y^0,Y)$ is
infinitely divisible with canonical measure $\rN$, standard computations
for Poisson measure yield that $\displaystyle
\rP_x(\tau_{Y^0}=m|\tau_Y=m)=\rN[\tau_{Y^0}=m|\tau_Y=m]$ that is 
\[
 \rP_x(\tau_{Y^0}=m|\tau_Y=m)=\exp{-\int_0^m
  \phi'(c(t))\; dt}.
\]
Notice that this formula  is also valid for the quadratic
branching mechanism ($\Pt(u)=\at  u+  \beta u^2$,
$\beta>0$, $\at\geq 0$), see Remark 5.3 in \cite{ad:cbpimcsbpi}. 

In  fact  this  formula  is  true in  a  general  framework.   Following
\cite{ad:cbpimcsbpi}, we consider the  branching mechanisms of the total
population and Eve-population are given by
\begin{align*}
   \Pt(\lambda)&=\at \lambda + \beta \lambda^2 +
\int_{(0,\infty )} 
\pt(d\ell) [\expp{-\lambda \ell} -1 + \lambda \ell],\\
\Pe(\lambda)&=\ae \lambda + \beta \lambda^2 +
\int_{(0,\infty )} 
\pe(d\ell) [\expp{-\lambda \ell} -1 + \lambda \ell],
\end{align*}
and the immigration function 
\[
\phi(\lambda)=\Pe(\lambda)-\Pt(\lambda)=
\aimm \lambda + \int_{(0,\infty )} \nu(d\ell) (1- \expp{-\lambda
  \ell}),
\]
where $\aimm=\ae-\at -\int_{(0,\infty )} \ell \nu(d\ell)\geq 0$ and
$\pt=\pe+\nu$, where $\pe$ and $\nu$ are Radon measures on $(0,\infty
)$ with $\int_{(0,\infty )} \ell \nu(d\ell)<\infty $. Notice the
condition $\int_{(0,\infty )} \ell \nu(d\ell)<\infty $ is stronger than
the usual condition on the immigration measure, $\int_{(0,\infty )}
(1\wedge\ell) \;  \nu(d\ell)<\infty $, but is implied by the requirement
that $\int_{(1,\infty )} \ell \nu(d\ell)<\int_{(1,\infty )} \ell
\pt(d\ell)<\infty $. 

Inspired by Theorem \ref{theo:rmax}, we consider 
$\displaystyle \cn(d\ell,dt,dz)= \sum_{i\in I} \delta_{(\ell_i, t_i,
  z_i)}(d\ell,dt,dz) $ a Poisson point measure with intensity  
\begin{equation}\label{eq:defcn}
\ind_{[0,m)}(t) \expp{-\ell c(m-t)}  \ell \left[\pe(d\ell)\delta_0(dz)+
  \nu(d\ell) \delta_1(dz) \right] dt.
\end{equation}
Intuitively,  the mark  $z_i$ indicates  if  the ancestor  (of the  last
individual alive)  alive at time $t_i$  had a new  mutation ($z_i=1$) or
not  ($z_i=0$). Note  however that  if $\beta>0$  we have  to  take into
account mutation on the  skeleton.  More precisely, let $T_1=\min \{t_i,
z_i=1\}$ be the first mutation on  the nodes in the ancestral lineage of
the last  individual alive  and let $T_2$  be an exponential  random time
with parameter  $\aimm$ independent of  $\cn$. The time $T_2$ corresponds
to  the first  mutation on
the skeleton for the ancestral  lineage of the last individual alive.  We
set  
\begin{equation}\label{eq:defT0}
\begin{cases}
T_0=\min(T_1,  T_2) & \mbox{ if }  \min(T_1,  T_2)<m,\\
T_0=+\infty  & \mbox{otherwise}.
\end{cases}
\end{equation}
In particular there is simultaneous extinction if and only if
$T_0=+\infty$.

For $t\ge 0$, let us denote by $\rN_t(dY^0,dY)$ the joint law of
$((Y^0(s-t),Y(s-t)), s\ge t)$ under $\rN$. 
Recall $\km$ given by \reff{eq:def-km}.
Conditionally on $\cn$ and $T_2$, let $\displaystyle  \sum_{j\in J} \delta_{(t_j, Y^{0,j}, Y^j)}$ be a Poisson point
measure,  with intensity 
\[
\km(dt) \rN_t[(dY^0,dY), \ \ind_{\{\tau_Y\le m\}}].
\]
We set 
\begin{equation}\label{eq:defy0'y'}
({Y'}^0,Y')=\sum_{t_j<T_0} (Y^{0,j}, Y^j)+  \sum_{t_j\geq T_0} (0, Y^j). 
\end{equation}
We write $\Q_m$ for the law  of $({Y'}^0,Y')$ computed for a given value
of $m$. 
\begin{theo}
   \label{th:Y0Ylaw}
   Under  $\Q_m$,  $({Y'}^0,Y')$   is  distributed  as  $(Y^0,Y)$  under
   $\rN[\cdot|\tau_Y=m]$  ,  or  equivalently,  under  $\int_0^{+\infty}
   |c'(m)|\Q_m(\cdot)dm  $, $({Y'}^0,Y')$  is  distributed as  $(Y^0,Y)$
   under $\rN$.
\end{theo}
Let us remark that this Theorem is very close to
Theorem \ref{theo:rmax} but only deals with CB and does not specify
the underlying genealogical structure. This is the purpose of a
forthcoming paper \cite{adv:pcrt} where the genealogy of multi-type CB
is described.

Intuitively, conditionally  on the last individual alive  being at time
$m$,  until the first  mutation in  the ancestral  lineage (that  is for
$t_j<T_0$) ,  its ancestors give birth  to a population  with initial Eve
type which has  to die before time $m$, and after  the first mutation on
the  ancestral  lineage  (that  is  for  $t_j\geq  T_0$),  there  is  no
Eve-population in  the descendants which  still have to  die before time
$m$.

Now, using that the distribution of $(Y^0,Y)$ is
infinitely divisible with canonical measure $\rN$, standard computations
for Poisson measure yield that $\displaystyle
\rP_x(\tau_{Y^0}=m|\tau_Y=m)=\rN[\tau_{Y^0}=m|\tau_Y=m]$. As 
\begin{align*}
\rN[\tau_{Y^0}=m|\tau_Y=m]
&=   \Q_m(T_0=+\infty )\\
&= \Q_m(T_1=+\infty )\Q_m(T_2\geq m)\\
&=\expp{-\int_0^m dt \int_{(0,\infty )} \expp{-\ell c(m-t)} \ell
  \nu(d\ell) } \expp{-\aimm m}\\
&=\expp{-\int_0^m dt \; \phi'(c(t))},
\end{align*}
we deduce the following Corollary. 
\begin{cor}[Probability of simultaneous extinction]
   \label{cor:SimDis}
We have for almost every $m>0$
\begin{equation*}
\rP_x(\tau_{Y^0}=m|\tau_Y=m)=\exp{-\int_0^m \phi'(c(t))\; dt}, 
\end{equation*}
where $c$ is the unique (non-negative) solution of \reff{eq:def-c}.
\end{cor}

\bigskip The paper is organized as follows. In Section \ref{sec:not}, we
recall some  facts on the  genealogy of the  CRT associated with  a Lévy
process. We  prove a Williams' decomposition for  the exploration process
associated  with the  CRT in  Section \ref{sec:Will}.  We  prove Theorem
\ref{th:Y0Ylaw}  in  Section  \ref{sec:gene}.  Notice  that  Proposition
\ref{prop:Ylaw} is a direct consequence of Theorem \ref{th:Y0Ylaw}.

\section{Notations}
\label{sec:not}

We  recall here  the construction  of the  Lévy continuum  random tree
(CRT) introduced in \cite{lglj:bplpep,lglj:bplplfss} and developed later
in \cite{dlg:rtlpsbp}. We  will emphasize on the height  process and the
exploration process which are the key tools to handle this tree. The
results of this section are mainly extracted from \cite{dlg:rtlpsbp}.

\subsection{The underlying Lévy process}\label{subsec:levy}
We consider  a $\R$-valued Lévy process $(X_t,t\ge
0)$ with  Laplace exponent $\psi$ (for  $\lambda\ge 0$
$\E\left[\expp{-\lambda X_t}\right]=\expp{t\psi(\lambda)}$) satisfying
\reff{eq:def_psi} and \reff{eq:contH}. 
Let $I=(I_t,t\ge 0)$
be  the infimum  process of  $X$, $I_t=\inf_{0\le  s\le t}X_s$,  and let
$S=(S_t,t\ge 0)$  be the supremum process,  $S_t=\sup_{0\le s\le t}X_s$.
We will  also consider for every $0\le  s\le t$ the infimum  of $X$ over
$[s,t]$:
\[
I_t^s=\inf_{s\le r\le t}X_r.
\]

The point 0 is regular for the Markov 
process $X-I$, and $-I$ is the  local time of $X-I$ at 0 (see
\cite{b:pl}, chap. VII). Let $\N$ be the associated excursion measure of
the process $X-I$ away from 0, and $\sigma=\inf\{t>0; X_t-I_t=0\}$ the length
of the excursion of $X-I$ under $\N$. We will assume that under $\N$,
$X_0=I_0=0$. 

Since $X$ is of infinite variation, 0 is also regular for the Markov 
process $S-X$. The local time, $L=(L_t, t\geq 0)$, of $S-X$ at 0 will be
normalized so that 
\[
\E[\expp{-\beta S_{L^{-1}_t}}]= \expp{- t \psi(\beta)/\beta},
\]
where $L^{-1}_t=\inf\{ s\geq 0; L_s\geq t\}$ (see also \cite{b:pl}
Theorem VII.4 (ii)).

\subsection{The height process and the Lévy CRT}
For each $t\geq 0$, we consider the reversed process at time $t$,
$\hat X^{(t)}=(\hat X^{(t)}_s,0\le s\le t)$ by:
\[
\hat X^{(t)}_s=
X_t-X_{(t-s)-} \quad \mbox{if}\quad  0\le s<t,
\]
and $\hat X^{(t)}_t=X_t$. The two processes $(\hat X^{(t)}_s,0\le s\le t)$
and $(X_s,0\le s\le t)$ have the same law. Let $\hat S^{(t)}$ be the
supremum process of $\hat X^{(t)}$ and $\hat L^ {(t)}$ be the
local time at $0$ of $\hat S^{(t)} - \hat X^{(t)}$ with the same
normalization as $L$. 

\begin{defi}[\cite{dlg:rtlpsbp}, Definition
  1.2.1 and Theorem 1.4.3]\label{def:height_process}
   There exists a process $H=(H_t, t\geq 0)$, called the
   height process, such that for all $t\geq 0$, a.s. $H_t=\hat
   L^{(t)}_t$, and $H_0=0$. Because of hypothesis \reff{eq:contH}, the
   height process  $H$ is continuous. 
\end{defi}

The  height process  $(H_t, t\in  [0,\sigma])$  under $\N$  codes  a
continuous  genealogical structure,  the Lévy  CRT, via  the following
procedure. 

\begin{itemize}
   \item[(i)] To each $t\in [0,\sigma]$ corresponds a vertex at
   generation $H_t$. 
   \item[(ii)] Vertex $t$ is an ancestor of vertex $t'$ if
   $H_t=H_{[t,t']}$, where 
\begin{equation}
   \label{eq:def_H}
H_{[t,t']}=\inf\{H_u, u\in [t\wedge
t', t\vee t']\}.
\end{equation}
In general $H_{[t,t']}$ is the generation of the last
   common ancestor to $t$ and $t'$. 
   \item[(iii)] We put $d(t,t')=H_t+H_{t'}- 2 H_{[t,t']}$ and identify $t$
   and $t'$ ($t\sim t'$) if $d(t,t')=0$. 
\end{itemize}
The Lévy CRT coded by $H$ is then the quotient set $[0,\sigma]/ \sim$,
equipped with the distance $d$ and the genealogical relation specified
in (ii). 

Let $(\tau_s, s\geq 0)$ be the right continuous inverse of $-I$:
$\tau_s=\inf \{ t>0; -I_t>s\}$. Recall that $-I$ is the local time of
$X-I$ at $0$. Let $L^a_t$ denote the local time at level $a$ of $H$
until time $t$, see
Section 1.3 in \cite{dlg:rtlpsbp}.

\begin{theo}[\cite{dlg:rtlpsbp}, Theorem 1.4.1]\label{theo:RK}
  The process $(L^a_{\tau_x}, a\geq 0)$ is under $\P$ (resp. $\N$)
  defined as $Y$ under $\rP_x$ (resp. $\rN$).
\end{theo}

In what follows, we will use  the notation $\N$ instead of $\rN$ for the
excursion measure to stress  that we consider the genealogical structure
of the branching process.

\subsection{The exploration process}
\label{sec:PLRT}
The height process is not Markov. But it is a simple
function of a measure-valued Markov process, the so-called exploration
process.

If $E$  is a Polish  space, let $\cb(E)$  (resp. $\cb_+(E)$) be the  set of
real-valued  measurable (resp. and non-negative) functions  defined on $E$
endowed    with   its   Borel    $\sigma$-field,   and    let   $\cm(E)$
(resp.  $\cm_f(E)$)  be  the  set  of  $\sigma$-finite  (resp.   finite)
measures  on  $E$, endowed  with  the  topology  of vague  (resp.  weak)
convergence.  For any measure $\mu\in\cm(E)$ and $f\in \cb_+(E)$, we write
$$\langle \mu,f\rangle =\int f(x)\,\mu(dx).$$

The     exploration    process     $\rho=(\rho_t,t\ge    0)$     is    a
$\cm_f(\R_+)$-valued process defined  as follows: for  every $f\in
\cb_+(\R_+) $,
$$\langle \rho_t,f\rangle =\int_{[0,t]} d_sI_t^sf(H_s),$$
or equivalently 
\begin{equation}\label{eq:def_rho}
\rho_t(dr)=\sum_{\stackrel{0<s\le t}
  {X_{s-}<I_t^s}}(I_t^s-X_{s-})\delta_{H_s}(dr)+\beta\ind_{[0,H_t]}(r)dr. 
\end{equation}
In particular, the total mass of $\rho_t$ is
$\langle \rho_t,1\rangle =X_t-I_t$.

For $\mu\in \cm(\R_+)$, we set 
\begin{equation}
   \label{def:H}
H(\mu)=\sup\, \Supp \mu,
\end{equation}
where $ \Supp \mu$ is the closed support of $\mu$,  with the
convention $H(0)=0$. We have

\begin{prop}[\cite{dlg:rtlpsbp}, Lemma 1.2.2]\label{prop:rho}
Almost surely, for every $t>0$,
\begin{itemize}
\item $H(\rho_t)=H_t$,
\item $\rho_t=0$ if and only if $H_t=0$,
\item if $\rho_t\neq 0$, then $\Supp \rho_t=[0,H_t]$. 
\end{itemize}
\end{prop}

In the definition  of the exploration process, as $X$  starts from 0, we
have  $\rho_0=0$ a.s.  To  state the  Markov property  of $\rho$,  we must
first define  the  process  $\rho$  started  at any  initial  measure  $\mu\in
\cm_f(\R_+)$.

For $a\in  [0, \langle  \mu,1\rangle ] $,  we define the  erased measure
$k_a\mu$ by
\[
k_a\mu([0,r])=\mu([0,r])\wedge (\langle \mu,1\rangle -a), \quad \text{for $r\geq 0$}.
\]
If $a> \langle  \mu,1\rangle $, we set $k_a\mu=0$.   In other words, the
measure $k_a\mu$ is the measure $\mu$ erased by a mass $a$ from the
top of
$[0,H(\mu)]$.

For $\nu,\mu \in \cm_f(\R_+)$, and $\mu$ with compact support, we
  define the concatenation $[\mu,\nu]\in \cm_f(\R_+) $ of the
  two measures by: 
\[
\bigl\langle [\mu,\nu],f\bigr\rangle =\bigl\langle \mu,f\bigr\rangle +\bigl\langle \nu,f(H(\mu)+\cdot)\bigr\rangle ,
\quad f\in \cb_+(\R_+).
\]

Finally,  we  set for  every  $\mu\in  \cm_f(\R_+)$  and every  $t>0$
$\rho_t^\mu=\bigl[k_{-I_t}\mu,\rho_t]$.  We say that $(\rho^\mu_t, t\geq
0)$  is  the  process  $\rho$  started at  $\rho_0^\mu=\mu$,  and  write
$\P_\mu$  for its  law. Unless  there is  an ambiguity,  we  shall write
$\rho_t$ for $\rho^\mu_t$.

\begin{prop}[\cite{dlg:rtlpsbp}, Proposition 1.2.3]
The process $(\rho_t,t\ge 0)$ is a càd-làg strong Markov process in
$\cm_f(\R_+)$.
\end{prop}

Notice that  $\N$ is  also the excursion  measure of the  process $\rho$
away  from $0$,  and  that $\sigma$,  the  length of  the excursion,  is
$\N$-a.e.  equal to $\inf\{ t>0; \rho_t=0\}$.

\subsection{The dual process and representation formula}
\label{sec:dual}

We  shall need the  $\cm_f(\R_+)$-valued process  $\eta=(\eta_t,t\ge 0)$
 defined by
\[
\eta_t(dr)=\sum_{\stackrel{0<s\le t}{X_{s-}<I_t^s}}(X_s-I_t^s)\delta
_{H_s}(dr)+\beta\ind_{[0,H_t]}(r)dr.
\]
The process $\eta$ is the dual process of $\rho$ under $\N$ thanks to
the following time reversal property: recall $\sigma$ denotes the
length of the excursion under $\N$.

\begin{prop}[\cite{dlg:rtlpsbp}, Corollary 3.1.6]\label{prop:timereverse}
The processes $((\rho_s,\eta_s); s\ge 0)$ and
$((\eta_{(\sigma-s)-},\rho_{(\sigma-s)-});$ $s\ge 0)$ have the same
distribution under $\N$.
\end{prop}

It also enjoys the snake
property: for all $t\geq 0, s\geq 0$
\[
(\rho_t,\eta_t)_{[0, H_{[t,s]})}=
(\rho_s,\eta_s)_{[0, H_{[t,s]})},
\]
that is the measures $\rho$ and $\eta$ between two instants coincide up
to the minimum of the height process between those two instants. 


We recall the Poisson representation of $(\rho,\eta)$ under $\N$. Let
$\mathcal{N}_*(dx\,   d\ell\,  du)$   be  a   Poisson  point   measure  on
$[0,+\infty)^3$ with intensity
$$dx\,\ell\pi(d\ell)\ind_{[0,1]}(u)du.$$
For every $a>0$, let us denote by $\mathbb{M}_a$ the law of the pair
$(\mu_a,\nu_a)$ of finite measures on $\R_+$ defined by:  for $f\in \cb_+(\R_+)$ 
\begin{align*}
\langle \mu_a,f\rangle  & =\int\mathcal{N}_*(dx\, d\ell\, du)\ind_{[0,a]}(x)u\ell f(x),\\
\langle \nu_a,f\rangle  & =\int\mathcal{N}_*(dx\, d\ell\, du)\ind_{[0,a]}(x)\ell(1-u)f(x).
\end{align*}
We finally set $\mathbb{M}=\int_0^{+\infty}da\, \expp{-\at 
  a}\mathbb{M}_a$. 

\begin{prop}[\cite{dlg:rtlpsbp},
  Proposition 3.1.3]
\label{prop:poisson_representation1}
For every non-negative measurable function $F$ on $\cm_f(\R_+)^2$,
\[
\N\left[\int_0^\sigma F(\rho_t, \eta_t) \; dt \right]=\int\mathbb{M}(d\mu\,
    d\nu)F (\mu, \nu) ,
\]
where $\sigma=\inf\{s>0; \rho_s=0\}$ denotes the length of the
    excursion. 
\end{prop}
We can then deduce the following Proposition. 
\begin{prop}\label{prop:poisson_representation3}
For every non-negative measurable function $F$ on $\cm_f(\R_+)^2$,
\[
\N\left[\int_0^\sigma F(\rho_t, \eta_t) \; dL^a_t \right]=\expp{-\at 
  a}\int \mathbb{M}_a(d\mu\,
    d\nu)F (\mu, \nu) ,
\]
where $\sigma=\inf\{s>0; \rho_s=0\}$ denotes the length of the
    excursion. 
\end{prop}

\section{Williams' decomposition}
\label{sec:Will}
We  work  under  the  excursion   measure.  As  the  height  process  is
continuous,  its supremum  $\hm=\sup\{  H_r; r\in  [0,\sigma]\}$ is attained.  Let
$\tm=\inf\{s\geq 0; H_s=\hm\}$. 

For every  $m>0$, we set $T_m(\rho)=\inf\{s>0,  H_s(\rho)=m\}$ the first
hitting time  of $m$ for  the height process.  When there is no  need to
stress the dependence  in $\rho$, we shall write  $T_m$ for $T_m(\rho)$.
Recall the function $c$ defined by \reff{eq:def-c} is equal to
\begin{equation}
   \label{eq:defcN}
c(m)=\N[T_m<\infty]=\N[\hm\le m].
\end{equation}

We set $\rho_d=(\rho_{\tm+s}, s\geq 0)$ and
$\rho_g=(\rho_{(\tm-s)+}, s\geq 0)$, where $x_+=\max(x,0)$.

For every finite measure with compact support $\mu$, we write
$\P_\mu^*$ for the law of the exploration process $\rho$ starting at
$\mu$ and killed when it first reaches 0.
We also set
$$\hat\P_\mu^*:=\lim_{\varepsilon\to 0}\P_\mu^*(\,\cdot\,|\, H(\mu)\le
\hm\le H(\mu)+\varepsilon).$$
We now describe the probability measure $\hat\P_\mu^*$  via a
Poisson decomposition. Let $(\alpha_i,\beta_i)$, $i\in I$ be the
excursion intervals of
the process $X-I$ away from 0 (well defined under $\P_\mu^*$ or under $\hat\P_\mu^*$). For every $i\in
I$, we define $h_i=H_{\alpha_i}$ and the measure-valued process $\rho^i$ by the formula
$$\langle
\rho^i_t,f\rangle=\int_{(h_i,+\infty)}f(x-h_i)\rho_{(\alpha_i+t)\wedge\beta_i}(dx).$$
We then have the following result.

\begin{lem}
\label{lem:Phat}
Under the probability $\hat\P_\mu^*$, the point measure $\displaystyle
\sum_{i\in I}\delta _{(h_i,\rho^i)}$ is a Poisson point measure with
intensity $\mu(dr)\N[\cdot,\, \hm\le m-r]$.
\end{lem}

\begin{proof}
We know (cf Lemma 4.2.4 of \cite{dlg:rtlpsbp}) that the point measure $\displaystyle
\sum_{i\in I}\delta _{(h_i,\rho^i)}$ is under $\P_\mu^*$ a Poisson point measure with
intensity $\mu(dr)\N(d\rho)$. The result follows then easily from
standard results on Poisson point measures.
\end{proof}

\begin{rem}\label{rem:approx_below}
Lemma \ref{lem:Phat}  gives also that, for every finite measure with compact
support $\mu$, if we write $\mu_a=\mu(\cdot\cap [0,a])$,
$$\hat\P_\mu^*=\lim_{a\to H(\mu)}\P_{\mu_a}^*(\,\cdot \,|\, \hm\le H(\mu)).$$
\end{rem}

\begin{theo}[Williams' Decomposition]
$ $
 \label{theo:rmax}
\begin{itemize}
\item[(i)] The  law of $\hm$  is characterized by  $\N[\hm\leq m]=c(m)$,
  where $c$ is the unique non-negative solution of \reff{eq:def-c}.
   \item[(ii)] Conditionally on $\hm=m$, the law of $(\Rm,\Em)$ is under $\N$
the law of
\[
\left(\sum_{i\in
  I}v_ir_i\delta_{t_i}+ \beta \ind_{[0,m]}(t) dt,\sum_{i\in
  I}(1-v_i)r_i\delta_{t_i}+ \beta \ind_{[0,m]}(t) dt\right),
\]
where $\sum\delta_{(v_i,r_i,t_i)}$ is a Poisson measure with intensity
$$\ind_{[0,1]}(v)\ind_{[0,m]}(t)\expp{-rc(m-t)}dv\, r\pi(dr)\, dt.$$
   \item[(iii)]   Under $\N$, conditionally on $\hm=m$, and $(\Rm,\Em)$, 
$(\rho_d, \rho_g)$ are independent and $\rho_d$ (resp. $\rho_g$) is
distributed as $\rho$ (resp. $\eta$) under  $\hat\P^*_{\Rm}$
(resp. $\hat\P^*_{\Em}$).    
\end{itemize}
\end{theo}

Notice (i) is a consequence of \reff{eq:defcN}. 
Point (ii)  is reminiscent of  Theorem 4.6.2 in  \cite{dlg:rtlpsbp} which
gives the description of the exploration process at a first hitting time
of the Lévy snake.

\bigskip The  end of this section is devoted  to the proof of  (ii) and
(iii) of this Theorem.

Let  $m>a>0$  be  fixed.  Let  $\varepsilon>0$.  Recall  $T_m=\inf\{t>0;
H_t=m\}$ is  the first hitting time  of $m$ for the  height process, and
set $L_m=\sup\{t<\sigma; H_t=m\}$ for the last hitting time of $m$, with
the    convention    that     $    \inf    \emptyset=+\infty    $    and
$\sup\emptyset=+\infty $.  We consider  the minimum of $H$ between $T_m$
and $L_m$: $\htl=\min\{ H_t; t\in [T_m, L_m]\}$.

We set $\rho^{(d)}=(\rho_{T_{\text{max},a}+t}, t\geq 0)$, with
$$T_{\text{max},a}=\inf\{t>\tm, H_s=a\},$$
  the path of
the exploration process on the right of $\tm$ after the hitting time of
$a$,  and $\rho^{(g)}=(\rho_{(L_{\text{max},a} -t)-}, t\geq 0)$, with 
$L_{\text{max},a}=\sup\{t<\tm; H_t=a\}$, the returned path of
the exploration process on the left of $\tm$ before its last hitting time of
$a$. Let us note that, by time reversal (see Proposition
  \ref{prop:timereverse}), the process $\rho^{(g)}$  is of the same type
  as $\eta$. This remark will be used later.

 To prove the Theorem, we shall compute 
\[
A_0=\N\left[F_1(\rho^{(g)}) F_2(\rho^{(d)}) F_3({\rho_\tm}_{|[0,a]})
  F_4({\eta_\tm}_{|[0,a]})   \ind_{\{ m\leq \hm<m+\varepsilon\}}
  \right]
\]
and  let  $\varepsilon$  go  down  to  $0$.  We  shall  see  in  Lemma
\ref{lem:approx}, that adding $ \ind_{\{\htl>a\}}$ in the integrand does
not change  the asymptotic  behavior as $\varepsilon$  goes down  to $0$.
Intuitively, if  the maximum  of the height  process is between  $m$ and
$m+\varepsilon$, outside a set of small measure, the height process does
not reach level $a$ between the  first and last hitting time of $m$.  So
that we shall compute first
\begin{equation}
   \label{eq:formuleAF}
A=\N\left[F_1(\rho^{(g)}) F_2(\rho^{(d)}) F_3({\rho_\tm}_{|[0,a]})
  F_4({\eta_\tm}_{|[0,a]})   \ind_{\{\htl>a, m\leq \hm<m+\varepsilon\}}
  \right].
\end{equation}

Notice that on $\{\htl>a\}$,  we have
$T_{\text{max},a}=T_{m,a}:=\inf\{s>T_m, H_s(\rho)=a\}$ and, from the snake
property,   ${\rho_\tm}_{|[0,a]}={\rho_{T_m}}_{|[0,a]}$ and 
${\eta_\tm}_{|[0,a]}={\eta_{T_m}}_{|[0,a]}$, so that 
\[
A=\N\left[F_1(\rho^{(g)}) F_2((\rho_{T_{m,a}+t}, t\geq 0))
  F_3({\rho_{T_m}}_{|[0,a]}) 
  F_4({\eta_{T_m}}_{|[0,a]})   \ind_{\{\htl>a, m\leq \hm<m+\varepsilon\}}
  \right].
\]
Let us remark that, we have
$$\ind_{\{\htl>a, m\leq \hm<m+\varepsilon\}}=\ind_{\{m\le
  \sup\{H_u,0\le u\le T_{m,a}\}<m+\varepsilon\}}\ind_{\{\sup\{H_u,u\ge
  T_{m,a}\}< m\}}.$$
By using the strong Markov property of the exploration process at time
$T_{m,a}$, we get 
\[
A\!\!=\!\!\N\!\left[F_1(\rho^{(g)})   F_4({\eta_{T_m}}_{|[0,a]}) \ind_{\{m\le
  \sup\{H_u,0\le u\le T_{m,a}\}<m+\varepsilon\}}
 F_3({\rho_{T_m}}_{|[0,a]})
  \E_{{\rho_{T_m}}_{|[0,a]}}^*\!\!\left[F_2(\rho)\ind_{\{\hm<
  m\}}\right]\right]
\]
and so, by conditioning, we get
\[
A=\N\left[F_1(\rho^{(g)})   F_4({\eta_{T_m}}_{|[0,a]}) 
  G_2({\rho_{T_m}}_{|[0,a]}) 
  \ind_{\{\htl>a, m\leq \hm<m+\varepsilon\}}
  \right],
\]
where  $G_2(\mu)=F_3(\mu)\E^*_\mu[F_2(\rho)|\hm<m]$. Using time
reversibility (see Proposition  \ref{prop:timereverse}) and the strong
Markov property at time $T_{m,a}$ again, we have
\begin{align*}
   A
&=\N\left[F_1(\rho^{(d)})   F_4({\rho_{T_m}}_{|[0,a]}) 
  G_2({\eta_{T_m}}_{|[0,a]}) 
  \ind_{\{\htl>a, m\leq \hm<m+\varepsilon\}}
  \right]\\
&=\N\left[G_1({\rho_{T_m}}_{|[0,a]}) 
  G_2({\eta_{T_m}}_{|[0,a]}) 
  \ind_{\{\htl>a, m\leq \hm<m+\varepsilon\}}
  \right],
\end{align*}
where $G_1(\mu)=F_4(\mu)\E^*_\mu[F_1(\rho)|\hm<m]$.

Now,  we use   ideas from  the  proof of  Theorem  4.6.2 of
\cite{dlg:rtlpsbp}.  Let us  recall the  excursion decomposition  of the
exploration process above level $a$.  We set
$\displaystyle \tau_s^a=\inf\left\{r,\ \int_0^rdu\,\ind_{\{H_u\le
    a\}}>s\right\}$. 
Let $\mathcal{E}_a$ be the $\sigma$-field generated by the process 
$(\tilde \rho_s,\ s\ge 0):=(\rho_{\tau_s^a},\ s\ge 0)$.
We also set $\tilde\eta_s=\eta_{\tau_s^a}$.

Let $(\alpha_i,\beta_i)$, $i\in I$ be the excursion intervals of $H$
above level $a$. For every $i\in I$ we define the measure-valued
process $\rho^i$ by setting
$$\begin{cases}
\langle \rho^i_s,\varphi\rangle =
\int_{(a,+\infty)}\rho_{\alpha_i+s}(dr)\varphi(r-a) & \mbox{if
}0<s<\beta_i-\alpha_i,\\
\rho_s=0 & \mbox{if }s=0\mbox{ or }s\ge \beta_i-\alpha_i,
\end{cases}$$
and the process $\eta^i$ similarly.
We also define the local time at the beginning of excursion $\rho^i$
by $\ell_i=L_{\alpha_i}^a$.
Then, under $\N$, conditionally on $\mathcal{E}_a$, the point measure
$$\sum_{i\in I}\delta_{(\ell_i,\rho^i,\eta^i)}$$
is a Poisson measure with intensity $\ind_{[0,L_\sigma^a]}(\ell)d\ell
\N[d\rho\ d\eta]$.

In particular, we have
\[
A=\N\left[\sum_{i\in I} \prod_{j\neq i} \ind_{\{T_m(\rho^j)=+\infty \}} 
G_1(\rho_{\alpha_i})
  G_2(\eta_{\alpha_i})
  \ind_{\{ m\leq \hm(\rho^i)<m+\varepsilon\}}
\right]. 
\]
Let us denote by $(\tau_\ell^a, \ell\geq 0)$ the right-continuous inverse of
$(L^a_s, s\geq 0)$. Palm
formula for Poisson point measures yields
\begin{align*}
A & =\N\left[\N\left[\sum_{i\in I} \prod_{j\neq i} \ind_{\{T_m(\rho^j)=+\infty \}} 
G_1(\rho_{\alpha_i})
  G_2(\eta_{\alpha_i})
  \ind_{\{ m\leq
  \hm(\rho^i)<m+\varepsilon\}}\biggm|\ce_a\right]\right]\\
&
  =\N\left[\int_0^{L_\sigma^a}d\ell\, G_1(\rho_{\tau^a_\ell})G_2(\eta_{\tau^a_\ell})\N[m\le
  \hm<m+\varepsilon]\N\left[\prod_{j\in
  I}\ind_{\{T_m(\rho^j)=+\infty\}}\Bigm|\ce_a\right]\right].
\end{align*}
A time-change then  gives 
\begin{equation}
   \label{eq:formuleAG}
   A
=v(m-a, \varepsilon) \N\left[\int_0^\sigma
  dL_s^a G_1(\rho_s)   G_2(\eta_s) \expp{- c(m-a) L^a_\sigma} \right], 
\end{equation}
where $v(x, \varepsilon)=c(x) - c(x+\varepsilon)=\N[ x\leq
\hm<x+\varepsilon]$. 
We have
\begin{align*}
   A
&=v(m-a, \varepsilon) \N\left[\int_0^\sigma
  dL_s^a G_1(\rho_s)   G_2(\eta_s) \expp{- c(m-a) L^a_s}\expp{- c(m-a)
    (L^a_\sigma -L^a_s)} \right]\\
&=v(m-a, \varepsilon) \N\left[\int_0^\sigma
  dL_s^a G_1(\rho_s)   G_2(\eta_s) \expp{- c(m-a) L^a_s}\expp{- \langle
    \rho_s, \N[1-\expp{ - c(m-a)L^{a- \cdot}_\sigma} ] \rangle}
\right],
\end{align*}
where we used for the last equality that the predictable projection of
$\expp{- \lambda (L^a_\sigma -L^a_s)}$ is given by 
$ \expp{- \langle
    \rho_s, \N[1-\expp{ - \lambda L^{a- \cdot}_\sigma} ] \rangle}$. 
Notice that by using the excursion decomposition above level $0<r<m$, we
have 
\[
c(m)=\N[T_m<\infty ]=\N[1-\expp{ - c(m-r)L^{r}_\sigma} ].
\]
In particular, we get 
\[
A=v(m-a, \varepsilon) \N\left[\int_0^\sigma
  dL_s^a G_1(\rho_s)   G_2(\eta_s) \expp{- c(m-a) L^a_s}\expp{- \langle
    \rho_s, c(m-\cdot) \rangle}\right].
\]
Using time reversibility, we have 
\[
A=v(m-a, \varepsilon) \N\left[\int_0^\sigma
  dL_s^a G_1(\eta_s)   G_2(\rho_s) \expp{- c(m-a) (L^a_\sigma-
    L^a_s)}\expp{- \langle 
    \eta_s, c(m-\cdot) \rangle}\right].
\]
Similar computations as those previously done give
\begin{align*}
   A
&=v(m-a, \varepsilon) \N\left[\int_0^\sigma
  dL_s^a G_1(\eta_s)   G_2(\rho_s) \expp{- \langle 
    \eta_s+\rho_s , c(m-\cdot) \rangle}\right]\\ 
&=v(m-a, \varepsilon) \N\left[\int_0^\sigma
  dL_s^a G_1(\rho_s)   G_2(\eta_s) \expp{- \langle 
    \rho_s+\eta_s , c(m-\cdot) \rangle}\right].
\end{align*}
Using Proposition \ref{prop:poisson_representation3}, we get
\[
A=v(m-a, \varepsilon) 
\expp{-\at 
  a}\int \mathbb{M}_a(d\mu\,
    d\nu) G_1(\mu)   G_2(\nu) \expp{- \langle 
    \mu+\nu , c(m-\cdot) \rangle}.
\]
We can give a first consequence of the previous computation.
\begin{lem}
\label{lem:approx}
   We have 
\[
\N[\htl>a,m\leq \hm<m+\varepsilon]=c'(m)\frac{c(m-a)
    -c(m-a+\varepsilon)}{c'(m-a)}.
\]
\end{lem}
\begin{proof}
Taking $F_1=F_2=F_3=F_4=1$ in \reff{eq:formuleAG}, we deduce
that 
\[
\N[\htl>a,m\leq
  \hm<m+\varepsilon]=v(m-a, \varepsilon) \N\left[
  L_\sigma^a  \expp{- c(m-a) L^a_\sigma} \right].
\]
Let    $a_0>0$   and   let    us   compute    $\displaystyle   B(a_0,a)=
\N\left[L_\sigma^a \expp{- c(a_0) L^a_\sigma} \right]$.
Thanks to Theorem \ref{theo:RK}, notice that 
\[
 B(a_0,a)=
\rN\left[Y_a \expp{- c(a_0) Y_a} \right]=\frac{\partial_{a_0} \rN[1-
  \expp{-c(a_0) Y_a}]}{c'(a_0)}.
\]
On the other hand, we have 
\[
c(a+a_0)=\rN[Y_{a+a_0}>0]=\rN[1-\rE_{Y_a}[Y_{a_0}=0]]=\rN\left[1-
\expp{- Y_a c(a_0)}\right],
\]
where we used  the Markov
property of $Y$ at time $a$ under $\rN$ for the second equality and 
\reff{eq:br-Y} with $\lambda$ going to infinity for the last. 
Thus, we get 
$\displaystyle B(a_0,a)=\frac{c'(a_0+a)}{c'(a_0)}$. We deduce that
\begin{align*}
   \N[\htl>a,m\leq
  \hm<m+\varepsilon]
&=v(m-a, \varepsilon) B(a-m,a)\\
&=c'(m)\frac{c(m-a)
    -c(m-a+\varepsilon)}{c'(m-a)}. 
\end{align*}
\end{proof}

Since $F_1, F_2, F_3$ and $F_4$ are bounded, say by $C$, we have 
 $|A-A_0| \leq  C ^4 \N[\htl<a,m\leq
  \hm<m+\varepsilon]$. From Lemma \ref{lem:approx}, we deduce that 
\[
\lim_{\varepsilon \rightarrow 0}\frac{|A-A_0|}{\N[m\leq
  \hm<m+\varepsilon]} \leq C^4 \left[1- \lim_{\varepsilon \rightarrow 0}
\frac{\N[\htl>a,m\leq
  \hm<m+\varepsilon]}{\N[m\leq
  \hm<m+\varepsilon]} 
\right] =0.
\]
We deduce that 
\begin{multline*}
   \lim_{\varepsilon \rightarrow 0} \frac{
\N\left[F_1(\rho^{(g)}) F_2(\rho^{(d)}) F_3({\rho_\tm}_{|[0,a]})
  F_4({\eta_\tm}_{|[0,a]})   \ind_{\{ m\leq \hm<m+\varepsilon\}}
  \right]}{\N[m\leq \hm<m+\varepsilon]} \\
\begin{aligned}
   &=\frac{c'(m-a)}{c'(m)} \expp{-\at 
  a}\int \mathbb{M}_a(d\mu\,
    d\nu) G_1(\mu)   G_2(\nu) \expp{- \langle 
    \mu+\nu , c(m-\cdot) \rangle}\\
&=\frac{ \int \mathbb{M}_a(d\mu\,
    d\nu) G_1(\mu)   G_2(\nu) \expp{- \langle 
    \mu+\nu , c(m-\cdot) \rangle}}{\int \mathbb{M}_a(d\mu\,
    d\nu)  \expp{- \langle 
    \mu+\nu , c(m-\cdot) \rangle}}\\
&= \int \tilde{\mathbb{M}}_a(d\mu\,
    d\nu) G_1(\mu)   G_2(\nu)\\
&= \int \tilde{\mathbb{M}}_a(d\mu\,
    d\nu) F_4(\nu) \E^*_\nu [F_1(\rho^{(d)})|\hm<m] 
F_3(\mu) \E^*_\mu [F_2(\rho^{(d)})|\hm<m] ,
\end{aligned}
\end{multline*}
where 
\begin{align*}
   \mu(dt)&=\sum_{i\in I} u_i \ell_i \delta_{t_i} +\beta
   \ind_{[0,a]}(t) dt\\
   \nu(dt)&=\sum_{i\in I} (1-u_i) \ell_i \delta_{t_i} +\beta
   \ind_{[0,a]}(t) dt,
\end{align*}
and $\sum_{i\in I} \delta_{(x_i,\ell_i,t_i)}$
 is under    $ \tilde{\mathbb{M}}_a$    a   Poisson  point   measure  on
$[0,+\infty)^3$ with intensity
$$\ind_{[0,a]}(t) dt\,\ell\expp{-\ell
  c(m-t)}\pi(d\ell)\ind_{[0,1]}(u)du.$$

Standard results on measure decomposition imply there exists a regular
version of the  probability  measure $\N[\,\cdot\,|\hm=m]$ and
that,  for almost every non-negative $m$,
\[
\N[\,\cdot\,|\hm=m]=\lim_{\varepsilon\to 0}\N[\,\cdot\,|m\le \hm
<m+\varepsilon].
\]
This gives (ii) and (iii) of Theorem \ref{theo:rmax} since $F_1,F_2,F_3,
F_4$ are arbitrary continuous functionals and by Remark
\ref{rem:approx_below}.

\section{Proof of Theorem \ref{th:Y0Ylaw}}
\label{sec:gene}

The  proof of  this Theorem  relies on  the computation  of  the Laplace
transform for $({Y'}^0,Y')$ and is given in the next three
paragraphs. The next paragraph gives some preliminary computations. 

\subsection{Preliminary computations}
\subsubsection{Law of $T_0$}\label{sec:law_T0}
Recall the definition of $\Q_m$ as the law of $({Y'}^0,Y')$
defined by \reff{eq:defy0'y'} and $T_0$ defined by \reff{eq:defT0} as
the first mutation undergone by the last individual alive.

For $r<m$, we have 
\begin{align*}
   \Q_m(T_0\in [r,r+dr], T_0=T_2)
&=\Q_m(T_2\in [r,r+dr])\Q_m(T_1>r)\\
&=dr\;\aimm \expp{-\aimm r} \exp{-\int_0^r dt \int_{(0,\infty )}
  \expp{-\ell c(m-t)}   \ell 
  \nu(d\ell)}\\
&=dr\;\aimm \expp{- \int_0^r \phi'(c(m-t))\; dt},
\end{align*}
and, with the notation  $\phi_0(\lambda)=\phi(\lambda)- \aimm \lambda$, 
\begin{align*}
   \Q_m(T_0\in [r,r+dr], T_0=T_1)
&=\Q_m(T_2>r)\Q_m(T_1\in [r,r+dr])\\
&=dr\; \phi'_0(c(m-r)) \expp{-\aimm r} \exp{-\int_0^r dt
  \int_{(0,\infty )} 
  \expp{-\ell c(m-t)}   \ell 
  \nu(d\ell)}\\
&=dr\;\phi'_0(c(m-r)) \expp{- \int_0^r \phi'(c(m-t))\; dt}.
\end{align*}
In particular, we have for $r<m$
\[
 \Q_m(T_0\in [r,r+dr])=dr\;\phi'(c(m-r)) \expp{- \int_0^r
   \phi'(c(m-t))\; dt}.
\]
and 
\begin{equation}
   \label{eq:QT>r}
 \Q_m(T_0>r)= \expp{- \int_0^r
   \phi'(c(m-t))\; dt}.
\end{equation}
Notice we have  $\displaystyle \Q_m(T_0=\infty )= \exp{- \int_0^m
   \phi'(c(t))\; dt}$. 
\subsubsection{Conditional law of $\cn$ given $T_0$}
\label{sec:condP}
Recall $\cn$ is under $\Q_m$ a Poisson point measure with intensity given by
\reff{eq:defcn}.
Conditionally on $\{T_0=r, T_0=T_2\}$, with $m>r>0$, $\cn$ is under
$\Q_m$  a point
Poisson measure with intensity
\begin{multline*}
   \ind_{[0,r)}(t) \expp{-\ell c(m-t)}  \ell \pe(d\ell)\delta_0(dz)dt +
 \\ 
 \ind_{(r,m)}(t) \expp{-\ell c(m-t)}  \ell 
\left[\pe(d\ell)\delta_0(dz)+
  \nu(d\ell) \delta_1(dz) \right] dt .
\end{multline*}

Conditionally on $\{T_0=r, T_0=T_1\}$, with $r<m$, $\cn$ is
distributed under
$\Q_m$ as $\tilde \cn+\delta_{(L,r,1)}$ where $\tilde\cn$ is a point
Poisson measure with intensity
\begin{multline*}
   \ind_{[0,r)}(t) \expp{-\ell c(m-t)}  \ell \pe(d\ell)\delta_0(dz)dt\\ 
+ \ind_{(r,m)}(t) \expp{-\ell c(m-t)}  \ell 
\left[\pe(d\ell)\delta_0(dz)+
  \nu(d\ell) \delta_1(dz) \right] dt,
\end{multline*}
and $L$ is a random variable independent of $\tilde\cn$ with distribution
$$\frac{\expp{-\ell c(m-r)} 
\ell \nu(d\ell)}{\int_{(0,\infty )}\expp{-\ell' c(m-r)}
\ell' \nu(d\ell')}.$$

Conditionally on $\{T_0=\infty \}$, $\cn$ is under $\Q_m$ a point
Poisson measure with intensity
\[
   \ind_{[0,m)}(t) \expp{-\ell c(m-t)}  \ell \pe(d\ell)\delta_0(dz)dt .
\]

\subsubsection{Formulas}
The following two formulas are straightforward: for all $x,\gamma\ge0$,
\begin{align}
\label{eq:intPt}
\Pe'(x+\gamma) -\Pe'(\gamma)&= 2\beta x+\int _{(0,\infty )} \expp{- \ell \gamma }\ell \pe(d\ell) [1- \expp{-
  \ell x}] ,\\
\Pt'(x+\gamma) -\Pt'(\gamma)&= 2\beta x+\int _{(0,\infty )} \expp{- \ell
  \gamma }\ell \pt(d\ell) [1- \expp{- 
  \ell x}] ,
\end{align}
Finally we deduce from \reff{eq:def-c} that $\Pt(c)=-c'$,
$\Pt'(c)c'=-c''$  and 
\begin{equation}
   \label{eq:edo-c}
\int \Pt'(c) = - \log(c').
\end{equation}

\subsubsection{Laplace transform}
\label{sec:Lap-trans}
 Recall $\tau_Y=\inf\{t>0; Y_t=0\}$ is the
extinction time of $Y$. Let  $\mue$ and $\mut$  be two finite  measures with
support  a  subset  of  a  finite  set  $A=\{a_1,  \ldots,  a_n\}$  with
$0=a_0<a_1<\cdots<a_n<a_{n+1}=\infty $.  For $m\in (0,+\infty )\setminus
A$, we consider
\begin{align*}
 w_m(t)&=\rN[1-\expp{-\int Y^0_{r-t} \;\mue(dr) - \int Y_{r-t}\;\mut(dr)
}\ind_{\{\tau_Y< m-t\}}],\\
w^*_m(t)&=\rN[1-\expp{- \int Y_{r-t}\;\mut(dr)
}\ind_{\{\tau_Y< m-t\}}].
\end{align*}
By      noticing       that      $\rN$-a.e.       $\ind_{\{\tau_Y<m-t\}}
=\lim_{\lambda\rightarrow\infty }  \exp{-\int Y_{r-t}\; \mu^\lambda(dr)}$,
where $\mu^\lambda(dr)=\lambda  \delta_m(dr)$, we deduce  from Lemma 3.1
in \cite{ad:cbpimcsbpi} that $(w_m,w^*_m)$  are right continuous and are
the unique non-negative solutions of : for $k\in\{0, \ldots, n\}$, $m\in
(a_k, a_{k+1})$, $t\in (-\infty , m)$,
\begin{equation}
 \label{eq:ww*1}
     w^*_m(t)+\int_{[t,a_k]} \Pt(w^*_m(r))dr 
= \int_{[t,a_k]} \mut(dr) +
c(m-a_k),
\end{equation}
\begin{multline}
\label{eq:ww*2}
   w_m(t)+\int_{[t,a_k]} \Pe(w_m(r))dr \\
= \int_{[t,a_k]} \mue(dr) +\int_{[t,a_k]}
\mut(dr) + 
c(m-a_k)+ \int_{[t,a_k]}  \phi(w_m^*(r))dr.
\end{multline}
We define 
\begin{equation}
   \label{eq:def-am}
\bar a_m=\max\{a_k; a_k<m, k\in \{0, \ldots,
n\}\}.
\end{equation}
Notice that 
$w_m(t)=w_m^*(t)=c(m-t)$ for $t\in (\bar a_m,m)$.

\subsection{Proof of Theorem \ref{th:Y0Ylaw}}
\subsubsection{Aim}
Theorem \ref{th:Y0Ylaw} will be proved as soon as we check that the
following equality
\[
w(0)=\int_{0}^\infty -c'(m) \Q_m[1-\expp{-\int {Y'}^0_{r}\;
  \mue(dr) - \int 
  Y'_{r} \;\mut(dr) 
}]dm 
\]
holds for all the possible choices of measures $\mue$ and $\mut$
satisfying the assumptions of Section \ref{sec:Lap-trans}, with
$w=w_\infty $ defined by  \reff{eq:ww*2}.

Notice the integrand of the right-hand side is null for $m<a_1$. Let
$\Delta$ denote the right-hand side. We have for $0<\varepsilon\leq a_1$:
\begin{align*}
   \Delta
&=\int_{\varepsilon}^\infty dm\;(-c'(m)) \Q_m[1-\expp{-\int {Y'}^0_{r}
  \;\mue(dr) - \int 
  Y'_{r}\;\mut(dr) 
}] \\
&= c(\varepsilon) + \int_\varepsilon^\infty dm\; \ind_{A^c}(m) c'(m)
\Q_m[Z],
\end{align*}
with, thanks to the definition \reff{eq:def-km} of $\km$, 
\[
Z=\exp{ -\int_0^{\bar a_m} \km(dt)\; [n_t 
\ind_{\{t<T_0\}} +n^*_t \ind_{\{t\geq T_0\}}]}
\]
and 
\begin{align*}
   n_t
&=\rN[( 1-\expp{-\int Y^0_{r-t} \; \mue(dr) -
    \int Y_{r-t}\;\mut(dr) 
} ) \ind_{\{\tau_Y\leq m-t\}}] 
=w_m(t) - c(m-t)\\
   n^*_t
&=\rN[(1-\expp{ -
    \int Y_{r-t}\;\mut(dr) 
} ) \ind_{\{\tau_Y\leq m-t\}}] =w^*_m(t) - c(m-t),
\end{align*}
with  $(w_m,w_m^*)$ the non-negative solutions of   \reff{eq:ww*1} and  \reff{eq:ww*2}.
Notice that   $w_m(t)=w_m^*(t)=c(m-t)$ for   $t\in (\bar a_m,m)$  and thus 
$n_t=n_t^*=0$ when  $t\in (\bar a_m,m)$.

We set $\displaystyle \Delta= c(\varepsilon)+\int_\varepsilon^\infty
\ind_{A^c}(m) (\Delta_1+\Delta_2+\Delta_3) \; dm$ with 
\begin{align*}
   \Delta_1&=c'(m) \Q_m[Z|T_0>\bar a_m]\Q_m(T_0>\bar a_m),\\
\Delta_2&=c'(m) \int_0^{\bar a_m} \Q_m[Z|T_0=r,T_0=T_1]  \Q_m(T_0\in [r,r+dr],
T_0=T_1),\\ 
\Delta_3&=c'(m) \int_0^{\bar a_m} \Q_m[Z|T_0=r,T_0=T_2]  \Q_m(T_0\in [r,r+dr],
T_0=T_2).
\end{align*}
We shall assume $m\not \in A$. 
\subsubsection{Computation of $\Delta_1$}
We have, using formula \reff{eq:def-km},
\begin{align*}
   \Delta_1
&=c'(m) \Q_m(T_0>\bar a_m)\Q_m[\expp{ -\int_0^{\bar a_m} \km(dt)
  n_t}|T_0>\bar a_m ] \\ 
&=c'(m) \expp{-\int_0^{\bar a_m}  \phi'(c(m-t))dt} \\
&\quad \exp\left\{-
 2 \beta \int_0^{\bar a_m} \!\!\!\!(w_m(t) -c(m-t))\; dt -    \int_0^{\bar a_m} \!\!\!\! dt \expp{-\ell
      c(m-t)} \ell \pe(d\ell) [1- \expp{- \ell (w_m(t) -
      c(m-t))}]\right\}\\
&=c'(m) \expp{-\int_0^{\bar a_m}\phi'(c(m-t))dt} \exp\left\{-
  \int_0^{\bar a_m}  dt [\Pe'(w_m(t))
  - \Pe'(c(m-t))]\right\}\\
&=c'(m) \expp{\int_{m-\bar a_m} ^m \Pt'(c(t))dt} \expp{- \int_0^{\bar a_m}
 dt\; \Pe'(w_m(t))}\\
&=c'(m-\bar a_m)  \expp{- \int_0^{\bar a_m} dt\; \Pe'(w_m(t))},
\end{align*}
where we used \reff{eq:edo-c} for the last equality to get 
\begin{equation}
   \label{eq:ePtc}
\expp{\int_{m-\bar a_m} ^m \Pt'(c(t))dt} =\expp{-
  [\log(c'(t))]_{m-\bar a_m}^m}=\frac{c'(m-\bar a_m)}{c'(m)}. 
\end{equation} 

\subsubsection{Computation of $\Delta_2$}
 Using Section \ref{sec:condP}, we get 
\begin{align*}
   \Q_m[Z|T_0=r,T_0=T_1] 
&=
\expp{-
 2 \beta \int_0^r (w_m(t) -c(m-t))\; dt-
 2 \beta \int_r^{\bar a_m} (w^*_m(t) -c(m-t))\; dt}\\
&\hspace{3cm} \exp( -\int_0^r dt \expp{-\ell c(m-t)} \ell \pe(d\ell) [1-\expp{-\ell
  n_t}])\\
&\hspace{3cm} \exp(
-\int_r^{\bar a_m }  dt \expp{-\ell c(m-t)} \ell \pt(d\ell) [1-\expp{-\ell
  n^*_t}]  )\\
&\hspace{3cm} \int_{(0,\infty )} \nu(d\ell') \frac{\expp{-\ell' c(m-r)}
  \ell'\expp{-\ell' n^*_r}}{\phi'_0(c(m-r))}\\
&=\exp(-\int_0^r dt [\Pe'(w_m(t)) - \Pe'(c(m-t))])\\
&\hspace{3cm} \exp(
-\int_r^{\bar a_m } dt [\Pt'(w_m^*(t))- \Pt'(c(m-t))  ])\\
&\hspace{3cm} \frac{\phi_0'(w^*_m(r))}{\phi_0'(c(m-r))}.
\end{align*}
We deduce from Section \ref{sec:law_T0}
\begin{align*}
\Delta_2
&=c'(m) \int_0^{\bar a_m} \Q_m[Z|T_0=r,T_0=T_1]  \Q_m(T_0\in [r,r+dr],
T_0=T_1),\\ 
&= c'(m) \int_0^{\bar a_m } dr\; \phi'_0(w^*_m(r)) \expp{- \int_0^r
  \phi'(c(m-t))\; dt} \exp(-\int_0^r dt [\Pe'(w_m(t)) - \Pe'(c(m-t))])\\
&\hspace{3cm} \exp(
-\int_r^{\bar a_m } dt [\Pt'(w_m^*(t))- \Pt'(c(m-t))  ])\\
&= c'(m)   \expp{ \int_0^{\bar a_m}
  \Pt'(c(m-t))\; dt} \int_0^{\bar a_m } dr\;  \phi'_0(w^*_m(r)) \expp{-\int_0^r dt\;
  \Pe'(w_m(t)) 
  - \int_r^{\bar a_m} dt\; \Pt'(w_m^*(t))}\\
&= c'(m-\bar a_m)   \int_0^{\bar a_m } dr\; \phi'_0(w^*_m(r))  \expp{-\int_0^r dt\;
  \Pe'(w_m(t)) -\int_r^{\bar a_m} dt\;
   \Pt'(w_m^*(t))}, 
\end{align*}
where we used \reff{eq:ePtc} for the last equality. 

\subsubsection{Computation of $\Delta_3$}
 Using Section \ref{sec:condP}, we get 
\begin{align*}
   \Q_m[Z|T_0=r,T_0=T_2] 
&=
\expp{-
 2 \beta \int_0^r (w_m(t) -c(m-t))\; dt-
 2 \beta \int_r^{\bar a_m} (w^*_m(t) -c(m-t))\; dt}\\
&\hspace{3cm} 
\exp\left\{ -\int_0^r dt \expp{-\ell c(m-t)} \ell \pe(d\ell) [1-\expp{-\ell
  n_t}] \right\}\\
&\hspace{3cm} 
\exp\left\{-\int_r^{\bar a_m }  dt \expp{-\ell c(m-t)} \ell \pt(d\ell)
  [1-\expp{-\ell   n^*_t}] \right\} \\
&=\exp(-\int_0^r dt [\Pe'(w_m(t)) - \Pe'(c(m-t))])\\
&\hspace{3cm} \exp(
-\int_r^{\bar a_m } dt [\Pt'(w_m^*(t))- \Pt'(c(m-t))  ]).
\end{align*}
We deduce from Section \ref{sec:law_T0} 
\begin{align*}
\Delta_3
&=c'(m) \int_0^{\bar a_m} \Q_m[Z|T_0=r,T_0=T_2]  \Q_m(T_0\in [r,r+dr],
T_0=T_2),\\ 
&= c'(m) \int_0^{\bar a_m } dr\; \aimm \expp{- \int_0^r
  \phi'(c(m-t))\; dt} \exp(-\int_0^r dt [\Pe'(w_m(t)) - \Pe'(c(m-t))])\\
&\hspace{3cm} \exp(
-\int_r^{\bar a_m } dt [\Pt'(w_m^*(t))- \Pt'(c(m-t))  ])\\
&= c'(m)   \expp{ \int_0^{\bar a_m}
  \Pt'(c(m-t))\; dt} \int_0^{\bar a_m } dr\; \aimm \expp{-\int_0^r dt\;
  \Pe'(w_m(t)) 
  - \int_r^{\bar a_m} dt\; \Pt'(w_m^*(t))}\\
&= c'(m-\bar a_m)   \int_0^{\bar a_m } dr\; \aimm \expp{-\int_0^r dt\;
  \Pe'(w_m(t)) -\int_r^{\bar a_m} dt\;
 \Pt'(w_m^*(t))}, 
\end{align*}
where we used \reff{eq:ePtc} for the last equality. 

\subsubsection{Computation of $\Delta_2+\Delta_3$}
We have
\[
   \Delta_2+\Delta_3 =c'(m-\bar a_m )  \int_0^{\bar a_m } dr
\;\phi'(w^*_m(r)) \expp{ -\int_0^r dt\; 
\Pe'(w_m(t))-\int_r^{\bar a_m } dt \;\Pt'(w_m^*(t))}.
\]

Differentiating  \reff{eq:ww*1}  w.r.t. time and $m$,
we get for  $t<m$
\[
\partial_m (w^*_m)'(t)- \partial_m w^*_m(t) \Pt'(w^*_m(t))=0.
\]
Notice also that for  $m>t\geq \bar a_m $, we have  $\partial_m
w^*(t)=c'(m-t)$ and thus 
\[
\partial_m
w^*(\bar a_m)=c'(m-\bar a_m).
\]
We get 
\[
\exp(-\int_r^{\bar a_m } dt \Pt'(w_m^*(t)))=\frac{\partial_m
  w^*_m(r)}{\partial_m w^*_m(\bar a_m )}= 
\frac{\partial_m
  w^*_m(r)}{c'(m-\bar a_m )}.
\]
Differentiating \reff{eq:ww*2}  
w.r.t. time and $m$, we get for  $t<m$ 
\[
\partial_m w_m'(t)- \partial_m w_m(t) \Pe'(w_m(t))=-\partial_m
w_m^*(t) \phi'(w^*_m(t)).
\]
We deduce that
\begin{align*}
    \Delta_2+\Delta_3
&= \int_0^{\bar a_m } dr
\partial_m
w_m^*(t)  \phi'(w^*_m(r)) \expp{ -\int_0^r dt\;
\Pe'(w_m(t))}\\
&=- \int_0^{\bar a_m } dr
[\partial_m w_m'(r)- \partial_m w_m(r)\Pe'(w_m(r))] \expp{ -\int_0^r dt\;
\Pe'(w_m(t))}\\
&= -\left[ 
\partial_m w_m(r) \expp{ -\int_0^r dt\;
\Pe'(w_m(t))}
\right]_0^{\bar a_m }\\
&= 
\partial_m w_m(0)- \partial_m w_m(\bar a_m ) \expp{
  -\int_0^{\bar a_m } dt \;
\Pe'(w_m(t))}.
\end{align*}

Notice also that for  $m>t\geq \bar a_m $ one has $\partial_m
w(t)=c'(m-t)$, in particular  $\partial_m w(\bar a_m)=c'(m-\bar a_m)$.
This implies that 
\[
\Delta_2+\Delta_3 =\partial_m w_m(0)- c'(m-\bar a_m) \expp{
   -\int_0^{\bar a_m } dt \;
 \Pe'(w_m(t))}.
\]

\subsection{Conclusion}
Thus, for $m\not\in A$,  we have 
\[
\Delta_1+\Delta_2+ \Delta_3=\partial_m w_m(0),
\]
and 
\[
  \Delta
= c(\varepsilon) + \int_\varepsilon^\infty \partial_m w_m(0)
=c(\varepsilon)+ w_\infty (0) - w_\varepsilon(0)
=w(0). 
\]
This ends the proof of the Theorem. 






\bigskip
{\bf Acknowledgments.} The authors wish to thank an anonymous referee
for his numerous and useful comments which 
improved  the  presentation of the paper.

\newcommand{\sortnoop}[1]{}

\end{document}